\documentclass[11pt]{article}

\usepackage[all]{xy}

\newtheorem{tm}{Theorem}[section]
\newtheorem{lm}[tm]{Lemma}
\newtheorem{pr}[tm]{Proposition}
\newtheorem{rmk}[tm]{Remark}
\newtheorem{cor}[tm]{Corollary}

\newtheorem{fact}[tm]{Fact}
\newtheorem{??}[tm]{Question}
\newtheorem{defi}[tm]{Definition}

\font\tenmsb=msbm10
\font\sevenmsb=msbm7
\font\fivemsb=msbm5
    
\newfam\msbfam
\textfont\msbfam=\tenmsb
\scriptfont\msbfam=\sevenmsb
\scriptscriptfont\msbfam=\fivemsb
\def\Bbb#1{{\fam\msbfam #1}}

\font\teneufm=eufm10
\font\seveneufm=eufm7
\font\fiveeufm=eufm5
\newfam\eufmfam
\textfont\eufmfam=\teneufm
\scriptfont\eufmfam=\seveneufm
\scriptscriptfont\eufmfam=\fiveeufm
\def\frak#1{{\fam\eufmfam\relax#1}}

\def\lorw{\longrightarrow}
\newcommand\n{\noindent}

\newcommand\ci{\cite}

\newcommand\rat{{\Bbb Q}}
\newcommand\comp{{\Bbb C}}

\newcommand\zed{{\Bbb Z}}

\newcommand\pn[1]{{\Bbb P}^{#1}}

\newcommand\blacksquare{{\hspace*{\fill} $\fbox{}$}}

\title{ A remark on singularities of primitive cohomology classes}
\author{
Mark Andrea A.  de Cataldo\thanks{
Partially supported by N.S.F.}\, 
and Luca Migliorini\thanks{ Partially supported by GNSAGA}
}
\date{November 2007}

\begin{document}\maketitle

\begin{abstract}
Green and Griffiths have introduced several notions of singularities
associated with normal functions, especially in connection with
middle dimensional primitive Hodge classes. 
In this note,  by 
using the more elementary aspects of the Decomposition Theorem,
we define  global and local
singularities associated with primitive middle dimensional
cohomology classes and  
by using the Relative Hard Lefschetz Theorem, we show
 that these singularities detect the global and local triviality of the
primitive class. In a final section, we write-up a  classical
inductive argument relating the Hodge conjecture to the local non-vanishing
of primitive classes.
\end{abstract}

\tableofcontents

\section{Introduction}
\label{intro}
In \ci{thomas}, Thomas has shown that the Hodge Conjecture (HC)
is equivalent to the existence of special hyperplane sections; see $\S$\ref{rtr}.
In \ci{grgr}, Green and Griffiths have introduced several notions of singularities
of normal functions, especially in connection with middle dimensional
primitive Hodge classes, and   have tied these notions to a positive
answer to  the HC. Roughly speaking, the HC is equivalent to the non vanishing
of these singularities at some point. As is well-known,  by virtue
of standard inductive arguments, the middle-dimensional case,
is the critical one.

The various   notions of singularities in \ci{grgr} reflect
the  subtle geometries which one needs to  explore to further this point of view. 
In this note, we 
concentrate on what is perhaps the most simple of these notions
and propose an essentially elementary  definition (cf. Definition \ref{definv})
of what we call the Green-Griffiths singularity  associated with a primitive
middle dimensional  class on a  projective manifold. 
There are a local (i.e. having to do with a hyperplane section)
and a global   version.
 Our  definition
is possible in view of the most basic properties
of the perverse filtration that can be read-off directly from
the Decomposition Theorem of Beilisnon, Bernstein, Deligne and Gabber
\ci{bbd}.
We prove, using the Relative Hard Lefschetz
Theorem \ci{bbd},  Propositions \ref{rhlsz} and \ref{rhlsz1} which 
establish that the local/global  triviality of the primitive
class is detected precisely by the Green-Griffiths invariants.

$\S$\ref{md} is for preliminaries on the Decomposition Theorem.
$\S$\ref{tggs} contains the results of this note mentioned above.
$\S$\ref{rhcsin} contains what we believe to be a classical
and well-known inductive argument concerning the HC. 
We could not find what we think is
a complete  reference for this result that puts into context
the notion of singularity introduced by Green and Griffiths.

We have reached the conclusions contained in $\S$\ref{tggs} in the Spring
of 2007, during our visit at I.A.S., Princeton. 

By using
 M. Saito's theories of  mixed Hodge modules and admissible normal functions,
 the   preprint  \ci{gg} introduces, among other things, a notion of singularity
that essentially  coincides with ours. Their methods and ours are entirely different
and we hope that the two different points of view can both be useful
for further geometric investigations.

We thank Phil Griffiths
for inspiring conversations.

\section{Decomposition Theorem formul\ae}
\label{md}
In this section we collect the  facts we need from the Decomposition
Theorem in \ci{bbd}. 

\subsection{Set-up}
\label{setup}
We work with rational cohomology. 
Let:

- $X^{2n} \subseteq  \pn{d}$ be a smooth, irreducible and  projective  manifold
of dimension $2n;$

- $\cal L$ be the hyperplane bundle; 
 
- $\zeta \in H^{2n}(X)$ be a  cohomology class; 

- there is the   universal hyperplane family (for simplicity, we write
${\Bbb P}$ also for ${\Bbb P}^{\vee}$)
$$
X \stackrel{q}\longleftarrow {\frak X}\stackrel{\pi}\longrightarrow \pn{d},
\qquad   \dim{{\frak X}} = {2n-1 +d}
$$
and we have the hyperplane sections as fibers
$$
\qquad
{\frak X}_p := \pi^{-1} (p), \quad
 \dim {\frak X}_p =2n-1.
$$

The  main object of investigation here  are   the classes 
$$
q^* \zeta \in H^{2n}({\frak X}), \qquad
(q^* \zeta)_{|   {\frak X}_p    }  =
 \zeta_{|   {\frak X}_p    }
 \in H^{2n}({\frak X}_p).
 $$
 
 Let $Z$ be a variety  $L$ be a local system on a dense
 open subset $U \subseteq Z_{reg}).$
The intersection cohomology complex $IC_Z(L)$ 
is a complex of sheaves on $Z.$
Its cohomology sheaves satisfy
\begin{equation}
\label{scsupp}
{\cal H}^l ( IC_Z(L)) = 0, \qquad l \neq [-\dim{Z}, -1].
\end{equation}
We have the intersection cohomology groups
\begin{equation}
\label{stic}
I\!H^k ( Z, L) := {\Bbb H}^{k -\dim{Z}} (Z, IC(L)).
\end{equation}
Clearly,
\begin{equation}
\label{obva}
I\!H^k (Z, L) =0, \qquad k \notin [0, 2\dim{Z}]. 
\end{equation}

 \subsection{The Decomposition Theorem  for $\pi$}
 \label{dt}
 The Decomposition Theorem  for  $\pi : {\frak X} \to \pn{d}$ gives a non canonical
 decomposition
 \begin{equation}
 \label{dtfac}
\phi: \;  \;
\bigoplus_{i \in \zed} \bigoplus_{j \in {\Bbb N}}
IC(L_{ij} ) [ -i - (2n-1+d)]   \; \simeq\;  R\pi_* \rat
\end{equation}
 where $L_{ij}$ is a local system on the codimension
 $j$ stratum $S_{d-j} \subseteq \pn{d}.$ Strata are not connected,
 so that by $IC$ here we mean the direct sum of the $IC$'s
 on the connected components of the same dimension.

The $L_{i0} $ are  the local systems on $S_d \subseteq  \pn{d} \setminus 
X^{\vee}:$
\begin{equation}
\label{dt3}
L_{i0} = R^{2n-1+i}: = (R^{2n-1+i}\pi_*\rat )_{| S_d   }.
\end{equation}

In what follows, the pedantic notation is to make the
formul\ae \, ready to use later.
We have  a non canonical decomposition for the  cohomology groups
\begin{equation}
\label{dt4}
H^{  l } ( {\frak X} ) \;\; = \;\;
\phi \; ( \;
\bigoplus_{i \in \zed}
\bigoplus_{j \in {\Bbb N} }
I\!H^{ (l-2n) -j - (i-1) } ( \overline{S_{d-j}}, L_{ij} ) \; ),
\end{equation}
one for cohomology sheaves
\begin{equation}
\label{dt5}
R^{l } \pi_* \rat  \;\;   = \;\; {\phi} \; (\;
\bigoplus_{i \in \zed }
\bigoplus_{j \in {\Bbb N} }
{\cal H}^{  (l-2n) -(d-j) -j - (i-1)} (   IC( L_{ij}  )) \; )
\end{equation}
and one for the cohomology groups
\begin{equation}
\label{dt6}
(   R^{ l   }    \pi_*  \rat) _p \; 
=\;  H^{  l } ( {\frak X}_p ) \;\;  = \; {\phi} \; (\; 
\bigoplus_{i \in \zed }
\bigoplus_{j \in {\Bbb N} }
{\cal H}^{  (l-2n) -(d-j) -j - (i-1)} (   IC( L_{ij}  ))_p \;)
\end{equation}

 \subsection{Definition of the filtrations on 
 $H({\frak X})$ and $H( {\frak X}_p)$}
 \label{fco}
 These filtrations  are discussed and used in our paper \ci{decmightam}.
 The theory of perverse sheaves filters the groups
$H( {\frak X}).$ 
The Decomposition Theorem  makes this  more visible.

\begin{rmk}
\label{remipf}
{\rm 
{\bf (Reminder on the perverse filtration)}
The perverse filtration on $H( {\frak X} )$
 is by Hodge substructures and 
it coincides, up to some shifts, with the monodromy filtration associated with
the nilpotent cup-product action of $q^*L$
on $H(\frak X)$  \ci{decmightam}. A splitting $\phi_{\cal L}$ in the category of 
Hodge structures exists, as shown in \ci{decmightadt},
but we do not need it here (see Remark  \ref{tho0}).
Over the regular part of $\pi$,  the perverse filtration coincides, up to re-numbering, 
with the filtration coming from the Leray spectral sequence.
}
\end{rmk}

The {\bf perverse filtration} on $H^{l} ( {\frak X})$
is the increasing filtration indexed by $i \in \zed:$ 
\begin{equation}
\label{smpfdef}
H^l_{\leq i} ( {\frak X} )\;  :=\; \; 
\phi \; (\; \bigoplus_{i' \leq i} \; \bigoplus_{j \in {\Bbb N} }
I\!H^{ (l-2n) -j - (i'-1) } ( \overline{S_{d-j}}, L_{i'j} ) \;) \;
\subseteq \; H^l ( {\frak X} ).
\end{equation}
The perverse filtration is  independent of the splitting $\phi.$
The graded pieces $H^l_i ( {\frak X})$  are canonically isomorphic to
\begin{equation}
\label{grpfex}
H^l_{ i} ( {\frak X} )\;  =\; \; 
  \bigoplus_{j \in {\Bbb N} }
I\!H^{ (l-2n) -j - (i-1) } ( \overline{S_{d-j}}, L_{ij} )  .
\end{equation}

The decomposition $\phi$   induces, in the same way an
increasing  filtration on the stalks which we call the 
{\bf  induced filtration} (also independent of $\phi$)
\begin{equation}
\label{infi00}
H_{\leq i}^{ l } ( {\frak X}_p) \; = 
 \; \phi \, (  \;
  \bigoplus_{i'  \leq i}  \;  \bigoplus_{j \in {\Bbb N}}
 {\cal H}^{   (l-2n) - (d-j) -j - ( i'-1)       }  (IC(L_{j} ) )  \;)_p
 \subseteq H^{ l }( {\frak X}_p),
 \end{equation}
 with graded pieces canonically isomorphic to 
 \begin{equation}
 \label{infi11}
 H_{ i}^{ l } ( {\frak X}_p)\; =  \;
 \bigoplus_{j \in {\Bbb N}}
 {\cal H}^{   (l-2n) - (d-j) -j - ( i-1)       }  (IC(L_{j} ) )  \; )_p.
 \end{equation}

 The restriction map $r: H({\frak X}) \to H({\frak X}_p)$ is filtered
 and strict (obvious since the map is a direct sum map)
 with respect to the perverse and induced filtrations:
 $$
r:  H_{\leq i}({\frak X}) \to H_{\leq i}({\frak X}_p).
 $$
 
 If we fix  a small neighborhood   $U \subseteq \pn{d}$ of $p $
 and set ${\frak X}_U: =\pi^{-1} (U),$ then we have 
 filtered isomorphisms
 \begin{equation}
 \label{filtiso}
r:  H_{\leq i}({\frak X}_U) \simeq H_{\leq i}({\frak X}_p).
\end{equation}
 
 \section{The Green-Griffiths singularity}
 \label{tggs}

 \subsection{Bound on the  filtrations on $H( {\frak X} )$
and $H( {\frak X}_p).$ }
 \label{bpf}
 In order to define the Green-Griffiths singularity we need
 (ii) below.
 
\begin{lm}
\label{bocgl}  {\rm (i)}
$$
H^{2n}_{\leq 1} ( {\frak X} ) =  H^{2n} ( {\frak X}), \qquad
H^{2n}_{\leq 1} ( {\frak X}_p ) =  H^{2n} ( {\frak X}_p ).
$$
{\rm (ii) }
If $\zeta \in H^{2n}(X)$ is primitive, then
$$
q^* \zeta \in  H^{2n}_{\leq 0} ( {\frak X} ),
\qquad
\zeta_{|  {\frak X}_p }     \in  H^{2n}_{\leq 0} ( {\frak X}_p ).
$$
\end{lm}
{\em Proof.}  By virtue of the obvious vanishing
(\ref{obva}) in negative degrees, the graded piece (\ref{grpfex})
\begin{equation}
\label{grpfexsn}
H^{2n}_{ i} ( {\frak X} )\;  =\; \; 
  \bigoplus_{j \in {\Bbb N} }
I\!H^{ -j - (i-1) } ( \overline{S_{d-j}}, L_{ij} ) \; = \; 0,
\qquad \forall \, i \geq 2.
\end{equation}
This holds also for  ${\frak X}_U.$
In view of (\ref{filtiso}),  this proves (i).
 
 \n
 The graded piece   
 $$
 H^{2n}_1( {\frak X} )  \; =\; 
 \bigoplus_{j\in {\Bbb N} } {I\!H}^{-j} ( \overline{S_{d-j}}, IC (L_{1,j}) )=
 I\!H^0 (\pn{d}, L_{10}).
$$
By (\ref{dt3}), $L_{10} =R^{2n}$ is the local system on the dense stratum
of $\pn{d}$ corresponding to the variation of $H^{2n}( {\frak X}_{\eta}),$
where ${\frak X}_{\eta}$ is a smooth hyperplane section.
The group in question is just the space of  global invariants.
 Since we are assuming that
$\zeta$ is primitive, $\zeta_{| {\frak X}_{\eta} }=0$ defines the zero section in this group
and (ii) follows.
\blacksquare

\subsection{Definition of the Green-Griffiths invariant}
\label{dggi}
We have the decompositions:
  (\ref{grpfex})
for $H({\frak X})$ (non canonical),
 (\ref{dt4})  for the graded 
 $H_i({\frak X})$ (canonical),
 (\ref{infi11}) for  
 $H({\frak X}_p)$ (non canonical) and (\ref{dt6}) for the graded
   $H_i({\frak X}_p)$ (canonical).
   
   Let $\zeta \in H^{2n}(X).$ There is  the  non canonical  decomposition 
   associated with $\phi:$
\begin{equation}
\label{dt33}
q^* \zeta \; =\;  {\phi} \, ( \;  \sum_{ij} [q^* \zeta]_{ij} \; ), \qquad
\zeta_{| {\frak X}_p } \; =\; \phi
\, (\; \sum_{ij} [  \zeta_{| {\frak X}_p }  ]_{ij} \; ),
\end{equation}
where the terms $[-]_{ij}$ depend on $\phi.$

Let $\zeta \in H^{2n}(X)$ be primitive. Then, by
Lemma \ref{bocgl},  the terms $[-]_0$ and
$[-]_{0j}$  are well-defined,
independently of $\phi:$
$$
[q^* \zeta]_0 \; = \; \sum_{j \in {\Bbb N}} 
 [ q^* \zeta ]_{0j} \; \in
H^{2n}_0 ( {\frak X} ) \; = \; 
\bigoplus_{j \in {\Bbb N} } {I\!H}^{-j+1} ( \overline{S_{d-j}}, IC(L_{0j})),
$$
$$
[ \zeta_{ | {\frak X}_p } ]_0  \; = \;  \sum_{j=0}^{d} [  
\zeta_{| {\frak X}_p}  ]_{0j}
\; \in H^{2n}_0 (  {\frak X}_p  )  \; = \;
\bigoplus_{j \in {\Bbb N} } 
{\cal H}^{-d+1  } (IC(L_{0j} ))_p.
$$

As in the proof of Lemma \ref{bocgl},  complemented
by the support condition (\ref{scsupp}), the terms
with  $j \geq 2$ are zero:
\begin{equation}
\label{dt444}
[q^* \zeta]_0 = [q^* \zeta]_{00} + [q^* \zeta]_{01}, \qquad
\end{equation}
\begin{equation}
\label{dt445}
[ \zeta_{ | {\frak X}_p } ]_0   =   
[  \zeta_{| {\frak X}_p}  ]_{00} +  
[  \zeta_{| {\frak X}_p}  ]_{01}
\end{equation}
where we write explicitly, remembering that our notation calls for
$IC_Z(L) $ to have cohomology sheaves in the interval $[-\dim{Z},-1]:$
\begin{equation}
\label{dt446}
[  \zeta_{| {\frak X}_p}  ]_{00} \in {\cal H}^{-d+1  } (IC( R^{2n-1} ))_p ,
\qquad
 [  \zeta_{| {\frak X}_p}  ]_{01} \in  {\cal H}^{-d+1  } (IC( L_{01} ))_p.
 \end{equation}

\begin{rmk}
\label{sga}
{\rm 
The local system $L_{01}$ is  usually defined on
the regular part of the dual variety $X^{\vee} \subseteq \pn{d}.$
If we take the embedding associated with $m {\cal L},$  $m \gg 0,$ then
the local system  $L_{01}=0.$ 
This follows from [SGA 7.2, XVIII.  5.3.5 (``Condition A'') and 6.4
(``Condition A" is verified for $m \gg 0)$]. In fact, the stalk  of $L_{01}$
at a general point of the dual hypersurface measures
the failure of the  adjunction map  5.3.2 (loc.cit.) to be an isomorphism
(it is surjective for Lefschetz pencils).
This  
can be also seen    by using the Clemens-Schmid sequence.
}
\end{rmk}

\begin{defi}
\label{definv}
Let $\zeta \in H^{2n}(X)$ be a primitive class.

\n
The \underline{global Green-Griffiths invariant $s(\zeta)$} is defined to be 
$$
s( \zeta) := [\zeta]_{00}  \in I\!H^1 ( \pn{d}, IC (R^{2n-1})).
$$
The \underline{local  Green-Griffiths invariant $s(\zeta)_p$} is defined to be 
$$
s( \zeta)_p := [\zeta_{| {\frak X}_p   }   ]_{00}  \in
{\cal H}^{-d+1} ( IC (R^{2n-1}))_p.
$$
\end{defi}

Clearly, these invariants depend on the embedding.

\begin{rmk}
\label{suppo}
{\rm From the conditions of support of $IC,$ it follows that the locus
$$
{\rm Sing } (\zeta ) \; := \; \{ p \in \pn{d}  \, | \;  s( \zeta)_p \neq 0  \}
$$
is of codimension at least two.
}
\end{rmk}

\begin{rmk}
\label{tho0}
{\rm
The Hodge-theoretic splitting    $\phi_{\cal L}$  of  \ci{decmightadt}
shows that $H^{2n}_0 (  {\frak X}   )$ is endowed with a natural 
pure Hodge structure  so  that if $\zeta$ is a Hodge class, then so are
$[\zeta]_0$ and  $s(\zeta).$
Our paper \ci{decmightadt} does not afford local results. However, 
using the M. Saito's general  theory of mixed Hodge modules 
one can reach similar conclusions for $[\zeta_{| {\frak X}_p}]_0$
and  $s(\zeta)_p$
}
\end{rmk}

\subsection{R. Thomas' result}
\label{rtr}
Thomas [2005] has proved that the Hodge conjecture is equivalent,
given an arbitrary middle dimensional  primitive Hodge class $\zeta,$
to the existence of $m \gg 0$ such that there exists  $p \in |m {\cal L}|$
with ${\frak X}_p$ nodal and $\zeta_{ |{\frak X}_p}\neq 0.$ 
If we drop the nodal requirement, the  resulting
weaker statement is of course still true.
Because of primitivity, the hypersurface must be singular
and nodality is an interesting improvement.

\subsection{The classes $[q^* \zeta]_0$
and $[\zeta_{| {\frak X}_p }]_0$ detect global/local triviality
}
\label{sprl}
Given a primitive class $\zeta \in H^{2n}(X),$ 
the class $q^* \zeta \in H^{2n}_{\leq 0}$ and it defines
a canonical element $[q^* \zeta]_0 \in H^{2n}_0 ( {\frak X}).$
Ditto for $\zeta_{| {\frak X}_p } \in H_{\leq 0}^{2n}
( {\frak X}_p )$ and 
$[\zeta_{| {\frak X}_p }]_0 \in
H_{ 0}^{2n}
( {\frak X}_p ) .$

The following proposition ensures that,
given {\em any} embedding $|{\cal L}|,$ the global/local
triviality of primitive classes
  is detected by the global/local classes  $[-]_0.$
  
It is a simple consequence of the Relative Hard Lefschetz Theorem
 \ci{bbd}.
In particular, we only need that cupping with 
the relatively ample line bundle
$q^* {\cal L}$ is injective on $H_{\leq -1}( {\frak X} )$ and on
$H_{\leq -1} (   {\frak X}_p )$ (see \ci{decmightadt}).

\begin{pr}
\label{rhlsz}
Let $\zeta \in H^{2n}(X)$ be primitive.

\n
(i) 
The class $\zeta =0$ IFF $[q^* \zeta]_0 =0.$

\n
(ii)
The class $\zeta_{|   {\frak X}_p  }  =0$ IFF $[\zeta_{
{\frak X}_p
}]_0 =0.$
\end{pr}
{\em Proof.} We prove (i). The proof for $(ii)$ is analogous.
One direction is trivial.
Let $\zeta \neq 0, $ so that $q^* \zeta \neq 0.$ 
If $[q^* \zeta]_0 =0,$ then $q^* \zeta \in H_{\leq -1}^{2n}
( {\frak X} ).$
By the Relative Hard Lefschetz, 
the cup product with $q^* {\cal L}$ is injective
on $H^{2n}_{\leq -1} ({\frak X}),$ contradicting
$q^* ({\cal L}  \cdot \zeta) =0.$
\blacksquare

\subsection{The local Green-Griffiths invariants detects
global/local triviality}
\label{lggidt}
 The Green-Griffiths invariant captures the primitive class, i.e.
we have the following, where by $m\gg 0,$ we mean that we replace
the embedding given by 
$|{\cal L}|$ with the one given by $| m {\cal L}  |, $ with $m \gg 0:$

\label{sprl1}
\begin{pr}
\label{rhlsz1}
Let $\zeta \in H^{2n}(X)$ be primitive.

\n
{\rm (i)}
The class $\zeta =0$ IFF $s(\zeta) =0.$

\n
{\rm (ii)} Let   $m\gg 0.$
The class $\zeta_{|   {\frak X}_p  }  =0$ IFF 
$s(\zeta)_p=0.$
\end{pr}
{\em Proof.} By Remark \ref{sga} and
by (\ref{dt444}) and (\ref{dt445}),  we 
have that, for $m \gg 0:$  
$$s(\zeta)_p = [ \zeta_{ | {\frak X}_p  } ]_{00} = [
\zeta_{|   {\frak X}_p   } ]_0
$$
and we apply Proposition \ref{rhlsz} to deduce (ii). The proof of
(i) is identical (except for the fact that the result holds
for every embedding).
\blacksquare

\begin{rmk}
\label{tho01}
{\rm
Proposition \ref{rhlsz1} implies that one can state R. Thomas'
result using 
$[\zeta_{|   {\frak X}_p }]_0,$ or 
$s(\zeta)_p,$ instead of
$\zeta_{|  {\frak X}_p}.$
We do not need the Hodge-theoretic nature of  $\phi_{\cal L}$
to prove the  resulting statement. One simply uses
$\zeta_{|  {\frak X}_p}$ and Deligne's mixed Hodge structures.
}
\end{rmk}

\begin{rmk}
\label{comp}
{\rm
If we understand correctly,  the content of $\S$6 of the paper \ci{gg}
by Brosnan et al. is equivalent to Remark \ref{tho01}.
}
\end{rmk}

\subsection{Further characterizations of vanishing of $s(\zeta)_{p}$}
\label{altre}
In this section we clarifiy  the relation between
$s(\zeta)_{p}$ and two other  invariants associated with $\zeta.$

The natural map $\rat_{ {\frak X}_p} \to IC_ {{\frak X}_p}$ induces
 a map $H^{2n}( {\frak X}_p) \to I\!H^{2n}( {\frak X}_p).$
By  \ci{decmightadt}, Thm. 3.2.1,  the kernel of this map
is precisely $W_{2n-1} H^{2n}( {\frak X}_p).$ Since $\zeta_{| {\frak X}_{p}}$ is of type $(n,n)$
for the mixed Hodge structure, 
it is not in the kernel and we 
have the following (cf.   \ci{grgr}, Thm. 2.ii  ): 

\begin{cor}
\label{icclass}
{\rm
Let  
$\zeta \in H^{n,n}_{{\rat}}(X)$  
be a primitive Hodge class on $X.$ 

\n
The class $ [ \zeta_{|   {\frak X}_{p}}]_{0}  =0  $
IFF $\zeta_{|  {\frak X}_{p}}=0$ 
 in
 $I\!H^{2n}( {\frak X}_p)= I\!H_{2n-2}( {\frak X}_p).$
 
 \n
Let $m \gg 0.$ The Green-Griffiths invariant  
 $s(\zeta)_p=0$
IFF  $\zeta_{|   {\frak X}_p  }  =0$
$I\!H^{2n}( {\frak X}_p)= I\!H_{2n-2}( {\frak X}_p).$
}
\end{cor}

It is possible to give a more precise characterization for the locus
in which the Green Griffiths invariant does not vanish, which appears as a natural generalization 
of the condition for the (non) extendability of a normal function in 
\ci{ezz}. We need  the preliminary

\begin{lm}
\label{hic}
{\rm     Let $U$ be a contractible neighborhood of $0 \in \comp^d,$ let $D\ni 0$ be a  divisor, and  
$L$ be a local system on $U \setminus D.$ Then 
${\cal H}^{1-d}(IC(L))_0$ injects naturally in ${\Bbb H}^1(U \setminus D, L).$}
\end{lm}
{\em Proof.}
Let $U^*=U_d \subseteq U_{d-1} \subseteq \ldots U_0=U$
be the ascending chain of open subsets  $U_{l}= \coprod_{{l'\geq l}}S_{l}$
associated with a stratification
of $(\comp^{d}, D),$ and denote by $j_{l}:U_{l+1} \to U_l$ the corresponding imbeddings.
We have the well-known formula
$$IC_U (L) \; :=\, \tau_{\leq -1} {Rj_0 }_{*} ( \ldots  ( \tau_{\leq -d+1} {Rj_{d-2}}_{*} 
( \tau_{\leq -d}  {Rj_{d-1}}_{*}
L [d] ) ) \ldots ).$$
Since the truncations relative to $j_i$ for $i\leq d-2$ are in degree bigger than or equal to $1-d,$ and   
$\tau_{\leq -d}  {Rj_{d-1}}_{*}
L [d] )=({j_{d-1}}_{*}
L) [d],$  setting $J:U_{d-1}\to U$ we have
$${\cal H}^{1-d}(IC(L))_0= {\Bbb H}^{1-d}(U, RJ_*{j_{d-1}}_{*}
L) [d]  )= H^1(U_{d-1},{j_{d-1}}_{*}L ).$$
The latter cohomology group is the term $E_2^{10}$ in the Grothendieck spctral sequence 
for ${\Bbb H}^1(U_{d-1},{Rj_{d-1}}_{*}L )=H^1(U^*,L),$ The statement
follows from the edge sequence.
\blacksquare

\medskip
Recall our notation ${\frak X}_{U}= \pi^{-1}(U).$
\begin{cor}
\label{punctured}
{\rm Let $U$ be a contractible neighborhood of $p,$ and 
$U^*=U \setminus U\cap X^{\vee}$. Let $m \gg 0.$
Then $s(\zeta)_{p}=0$
IFF  $q^*\zeta_{|  {\frak X}_{U^*}}  =0$ in 
$H^{2n}(   {\frak X}_{U^*}).$} 
\end{cor}
{\em Proof.}
Let $u:U^* \to U$ be the open imedding. The map $\pi: \pi^{-1}(U^*) \to U^{*}$
is smooth so that we have, by Deligne's theorem, the splitting on the right end side below:
 $$
R\pi_*\rat_{ {\frak X}_{U}} \lorw  Ru_*u^* R\pi_*\rat_{{\frak X}_{U}}
\simeq \oplus Ru_*R^l\pi_*\rat_{{\frak X}_{U^*}}[-l].
$$
The statement follows from  \ref{rhlsz1}.ii  
and \ref{hic} applied to $L= R^{{2n-1}}= R^{2n-1}\pi_*\rat_{{\frak X}_{U^*}}.$
\blacksquare

\section{Relation between HC and $s(\zeta)_p$ }
\label{rhcsin}
The content of this section is well-known and  is essentially a 
re-writing of  an argument of  R. Thomas to be found in \ci{thomas}.
We do not fully understand
Thomas's  Hilbert scheme argument: we are un-able to rule out
the existence of  parasite contributions in the cycle
he constructs that do not ensure that it
restricts to the hyperplane sections as wanted.

 Let us summarize  briefly what follows.
If the HC holds, after passing to $|m {\cal L}|,$ $m \gg 0,$
a primitive  Hodge class $ \zeta $ has 
$\zeta_{| {\frak X}_p } \neq 0$ for some point $p \in \pn{d_m}.$  
In fact, assuming HC, we find an effective cycle pairing non-trivially
with $\zeta.$ Then  we find a hypersurface $V \in |m {\cal L}|,$
$m \gg 0,$   containing that cycle.
Then $\zeta_{|V} \neq 0.$ This $V$ gives the wanted $p \in \pn{d_m}.$

Conversely, using an inductive argument,
if  $\zeta_{| {\frak X}_p } \neq 0$ for some $p \in \pn{d},$
then one produces an algebraic $n$-cycle
in ${\frak X}_p$ that pairs non trivially with $\zeta.$ This means
that the arbitrary, primitive Hodge class  $\zeta$  is not
perpendicular to algebraic cycles and this implies HC. 

What follows is presumably classic and well-known.

\subsection{HC and ${A^k}^{\perp}$}
\label{hcpperp}
Let $Y$ be a  smooth projective manifold of dimension $d.$
Fix $k \in \zed.$
We have 
$$
A^k(Y) \subseteq H^k_g(Y) \subseteq H^{k,k}(Y) \subseteq H^{2k}_{\comp}(Y),
$$
with $H^k_g(Y)$ the Hodge classes i.e.
the pure Hodge substructure  given by  $H^k_g (Y) = H^{k,k}(Y) \cap H^{2k}_{\rat}
(Y)
\subseteq H^{2k}_{\rat}(Y),$ and with $A^k(Y)$ the algebraic
classes, i.e. the image of the cycle class map.
By Hard Lefschetz, there are   isomorphisms
$$
H^k_g (Y) \simeq H^{d-k}_g(Y).
$$
There is the cup product  bilinear map
$$
H^k_g (Y)  \times H^{d-k}_g (Y)  \lorw \rat, \qquad (a,b) \mapsto \int_Y a\cup b
$$
This  bilinear map is nondegenerate.
 Assume  that
$2k \leq d,$  else $2(d-k) \leq d$ and we can switch $a$ with $b$ in what follows. 
Since $\dim_{\rat} H^k_g(Y)=\dim_{\rat} H^{d-k}_g(Y), $ it is enough to show
that if $a\neq 0,$ then $\int_Y a\cup -$ is not the zero map. This follows
from the primitive decomposition and  the Hodge Riemann bilinear relations.

We have the following immediate 
\begin{fact}
\label{akhk}
The Hodge conjecture HC for $Y$ is equivalent to an orthogonality statement:
$$
A^k(Y) \, =\,  H^k_g(Y)  \qquad \mbox{IFF} \qquad
 0= (A^k(Y))^{\perp}   \subseteq H^{d-k}_g (Y).
$$
\end{fact}

\begin{lm}
\label{hc=perp}
The following are equivalent:

\n
1) HC holds;

\n
2) for every smooth projective  even dimensional
$X^{2n}$,  we have $A^n(X)^{\perp} =0.$

\end{lm}
{\em Proof.}
1) $\Longrightarrow$ 2): it is Fact \ref{akhk}.

\n
2) $\Longrightarrow$ 1).  We fix $Y^d$  and $k \in \zed.$
We want to prove that $A^{d-k}(Y) = H^{d-k}_g(Y).$ 
By Fact \ref{akhk}, this is equivalent to showing   that
$A^{d-k}(Y)^{\perp} = 0.$ 
Let $ a \in A^{d-k}(Y)^{\perp}.$
 We 
need to show that
$a=0.$ 
If $d=2k,$ then we are done by assumption 2).
We have two cases: $d < 2k$ and $d > 2k.$
By contradiction, assume $a \neq 0.$

\n
Assume $d  < 2k.$  Take $X^{2k} :=Y^d \times T^{2k-d},$ $T$ any projective manifold of dimension  $2k-d$ and let  $p : X \to Y$ be the projection.
Note that $p^*$ is injective so that  $0 \neq p^* a \in H_g^k(X).$
Let $Z \in A^k(X).$ Since, by the projection formula and the choice of $a,$
we have that  $p^*a \cdot Z = a \cdot 
p_* Z=0,$ we see that $p^*a \in A^k(X)^{\perp}.$ By 2), $p^*a=0$ and
this is a contradiction.

\n
Assume $d >2k.$ Slice down $Y^d$ to $X^{2k}$ by taking
the complete intersection  of $d-2k$ general hyperplane sections. Weak Lefschetz implies that
$0 \neq a_{|X} \in H^k_g(X).$ By  2), we have that
there must be an algebraic cycle $Z$  on $X$ such that
$a_{|X} \cdot Z \neq 0.$ This implies that   $a \cdot Z \neq 0$ (the product is now on
$Y$).
This contradicts the assumption $a \in A^{d-k}(Y)^{\perp}.$
\blacksquare

\subsection{HC and $\zeta_{| {\frak X}_p} $}
\label{hcsz}
\begin{lm}
\label{hcszea}
Let $X^{2n}$ be smooth and projective, $\zeta \in H^{2n}_g.$
Assume that HC(X) holds.
Then there are  $m \gg 0,$ and $ p \in \pn{d_m} = |m {\cal L}|$
such that 
$$
 \zeta_{ | {\frak X}_p} \neq 0.
 $$
If, in addition, $\zeta$ is primitive, then 
 the Green-Griffiths invariant
$$
s(\zeta)_p \neq 0.
$$
\end{lm}
{\em Proof.}
The pairing $H^n_g(X) \times H^n_g(X)\to \rat$ is perfect.
It follows that there is $a \in H^n_g(X)$ such that
$\zeta \cdot a \neq 0.$ By HC, $a = \sum r_i Z_i,$ $Z_i$ cycle classes.
It follows that there is an irreducible codimension $n$
cycle $Z$ with $\zeta \cdot Z \neq 0.$ This means that $\zeta_{|Z} \in
H^{2n}(Z) = H_{2n}(Z)^{\vee}$ is non-zero.
By Serre, there is $m \gg 0$ such that there is a,
necessarily singular,  element ${\frak X}_p$
of $|m {\cal L}|$
containing $Z.$ 
 It follows that $\zeta_{ |{\frak X}_p} \neq 0.$
We conclude by Proposition  \ref{rhlsz1}.
\blacksquare

\begin{rmk}
\label{bana}
{\rm
If $\zeta$ is not primitive, then
$\zeta_{| {\frak X}_p} \neq 0$ for general $p.$
The point is to realize this for primitive $\zeta,$ for which
the above restriction is zero for general $p$ and every
 $m .$
}
\end{rmk}
\begin{pr}
\label{hcsg}
Assume that for every $X^{2n},$  $0 \neq \zeta \in H^n_g(X)$
there is $m\gg 0$ and $p \in |m {\cal L}|$ with $[\zeta_{| \frak{X}_p }] \neq 0.$
Then HC holds for projective manifolds.
\end{pr}
{\em Proof.} 
We want to prove HC($d$) for every $d.$
The case $d=0$ is trivial. Assume we have done the case $d-1.$

\n
By taking a smooth hyperplane section we deal with $H^k(X), $
$k \geq d+1.$ 

\n
Using a Lefschetz pencil 
we deal with $H^k(X),$ $k \leq d-1.$  Here is how.
Take a nice Lefschetz pencil (condition A of SGA is ok)
$ X \stackrel{u}\leftarrow\widetilde{X} \to \pn{1}.$ 
By the blow-up formula, HC holds for 
$X$ IFF it holds for $\widetilde{X}.$ So we work on $H^k(\widetilde{X}).$
Let $v:X_t
\to \widetilde{X}$ be the general member.
We get an  isomorphism  of  pure Hodge structures
\begin{equation}
\label{1}
 H^{k-2}(X_t)(-1) \oplus H^k(X) \stackrel{g_*,u^*}\simeq H^k(\widetilde{X}) 
\end{equation}
($H^j(X) = H^{j}(X_t)^{\pi_1}$).
By induction we have HC for $X_t$ and $g_*$ sends
cycle classes to cycle classes.   Let $a \in H^k_g(X).$ Then
$u^* a \in H^k_g(\widetilde{X})$ and $u^*a_{| X_t} \in H^k_g (X_t).$ 
By HC for each smooth $X_t,$ we have 
$$
N_1 u^*a_{| X_t} = N_3 Z'_t - N_2 Z''_t
$$
with $N_i \in {\Bbb N}^*$ and the $Z's$ are effective cycle classes
on $X_t.$   For $l \gg 0,$ we have that $Z''_t$ is contained in 
a hypersurface of $|l {\cal L}|.$ By iterating, we have that
$Z''_t$ is contained in the complete intersection
of $k/2$ hypersurfaces in $|l {\cal L}|,$ $l \gg 0.$ It follows that
$$
Z''_t = \nu {\cal L}^{k/2} - T_t
$$
with $\nu  \gg 0$ and $T_t$ effective on $X_t.$
It follows that
$$
N_1 u^* a_{| X_t} + N_2 \nu {\cal L}^{k/2} = N_3 Z'_t + N_2 T_t = {\cal Z}_t,
$$
where ${\cal Z}_t$ is an effective cycle (class).
The set of $t$'s is uncountable.  
Take $Hilb(\widetilde{X}/\pn{1}).$ 
It has a countable number of components. It follows that there is an irreducible component of it that maps onto $\pn{1}$ and that it contains
an uncountable number of the $Z_t$'s (which then share necessarily the
$N_1$ and $N_2$). This irreducible component contains a curve
$C$ mapping surjectively to $\pn{1}$ and passing through a point corresponding
to one of the ${\cal Z}_t.$ We take that  $t$ and that $X_t.$
Take the base change family over $C.$ 
This gives an algebraic cycle $\widetilde{\cal Z}'$
on $\widetilde{X}$ whose restriction to $X_t$
 ${\cal Z}_t'$ has  cohomology class  a positive
 $N$ multiple of ${\cal Z}_t.$ By (\ref{1}) 
 $$
 {\cal Z}'  =  g_*\alpha + u^* \beta, \quad 
 {\cal Z}'_{| X_t} = N {\cal Z}_t =   NN_1 u^* a_{| X_t} 
 +  N N_2 \nu u^* {\cal L}^{k/2}, \qquad 
 g_* \alpha_{|X_t} =0.
 $$
 Since $u^*$ is injective in this range by Weak Lefschetz, 
 $$
 \beta =  NN_1a + NN_2  \nu {\cal L}^{k/2}.
 $$
 So, on $\widetilde{X}:$
 $$
 {\cal Z}' = g_* \alpha   + u^* (N N_1 a + N N_2 \nu {\cal L}^{k/2}).
 $$
 The summand $g_* \alpha$ is algebraic by induction.
 It follows that $u^*a$ is algebraic and we are done:
 we have dealt using induction with $H^k_g(X),$ $k<d-1.$
 
 \n
 The argument above looks like the one in 
 \ci{thomas}, however it is different. 
 We felt the need to write this argument up, since
\ci{thomas}
 does not  seem to deal with the possible parasites extra components
 that may arise in that construction after one restricts to a hyperplane
 section of the pencil.

\n
If $d$ is odd, then we are done.
If $d$ is even, the remaining case is $H^{2n}(X^{2n}).$
We deal with this case.  By the  inductive hypothesis, we know that 
HC($2n-1$) is ok.
By Lemma \ref{hc=perp} it is enough  to show
that for every $X^{2n}$ we have that $A^n(X)^{\perp} =0,$
i.e.  that a non zero Hodge class $\zeta  \in H^n_g(X)$
 cannot be perpendicular to algebraic cycles.

\n
We have $m\gg 0$ and $p \in |m {\cal L}|$ with $\zeta_{| {\frak X}_p} \neq 0.$

\n
Take a desingularization $f: \widetilde{{\frak X}_p} \to {\frak X}_p.$
Note that the domain may be disconnected.
By mixed Hodge theory, $H^{2n}( {\frak X}_p) \to H^{2n}
( \widetilde{{\frak X}_p})$ has as kernel the classes of weights
$\leq 2n-1.$ In particular, $f^* \zeta \neq 0 \in H^n_g 
(\widetilde{ \frak{X}_p}).$

\n
By the inductive hypothesis HC(2n-1), 
the  HC holds on $\widetilde{{\frak X}_p}.$
By Fact \ref{akhk}, there is a cycle class 
$W \in A^{n-1}( \widetilde{{\frak X}_p} ) $ 
such that 
$ g^*\zeta \cdot W \neq 0,$  where 
$g: \widetilde{{\frak X}_p} \to X.$
It follows that $\zeta \cdot g_* W \neq 0,$ i.e. we found a cycle
class which is not perpendicular to $\zeta.$
\blacksquare

Authors' addresses:

\smallskip
\n
Mark Andrea A. de Cataldo,
Department of Mathematics,
SUNY at Stony Brook,
Stony Brook,  NY 11794, USA. \quad 
e-mail: {\em mde@math.sunysb.edu}

\smallskip
\n
Luca Migliorini,
Dipartimento di Matematica, Universit\`a di Bologna,
Piazza di Porta S. Donato 5,
40126 Bologna,  ITALY. \quad
e-mail: {\em migliori@dm.unibo.it}

\end{document}